\documentclass[leqno]{article}
\usepackage{amsmath,amsthm,amsfonts,amssymb}
\usepackage{mathrsfs}

\allowdisplaybreaks

\newfam\msbfam
\font\tenmsb=msbm10    \textfont\msbfam=\tenmsb \font\sevenmsb=msbm7
\scriptfont\msbfam=\sevenmsb \font\fivemsb=msbm5
\scriptscriptfont\msbfam=\fivemsb

\newfam\bigfam
\font\tenbig=msbm10 scaled \magstep2   \textfont\bigfam=\tenbig
\font\sevenbig=msbm7 scaled \magstep2 \scriptfont\bigfam=\sevenbig
\font\fivebig=msbm5 scaled \magstep2
\scriptscriptfont\bigfam=\fivebig

\begin{document}

\title{{Boundedness of Singular Integral Operators on Weak Herz Type Spaces with Variable Exponent}
\footnotetext{\noindent *Corresponding author.} }

\author{Hongbin WANG{$^{*}$} and Zongguang LIU}

\date{}

\maketitle

\begin{minipage}{13.5cm}
{\bf Abstract}\quad  \small{In this paper, the authors define the weak Herz spaces and the weak Herz-type Hardy spaces with variable exponent. As applications, the authors establish the boundedness for a large class of singular integral operators including some critical cases.}

\medskip


\medskip

{\bf MR(2010)Subject Classification}\quad 42B35; 42B20; 46E30.

\end{minipage}

\section{ Introduction\label{s1}}$\indent$

The main purpose of this article is to introduce the weak Herz type spaces with variable exponent and give the boundedness for a large class of singular integral operators. The classical weak Herz type spaces can be traced to the work of Hu, Lu and Yang ([7,8]) on the study of the boundedness for some operators. The classical weak Herz type spaces mainly include the weak Herz spaces and the weak Herz-type Hardy spaces. As well known, many important operators in harmonic analysis, such as Calder\'{o}n-Zygmund operators, Ricci-Stein oscillatory singular integral operators, are not bounded on $L^1(\mathbb{R}^{n})$, but they map $L^1(\mathbb{R}^{n})$ to weak $L^1(\mathbb{R}^{n})$. So the weak Lebesgue space $WL^p(\mathbb{R}^{n})$ are very important in operator theory. The weak Herz spaces are the analog of $WL^p(\mathbb{R}^{n})$ in the setting of Herz spaces. Meanwhile, the weak Herz-type Hardy spaces are the suitable local version of the weak Hardy spaces $WH^p(\mathbb{R}^{n})$. The classical weak Herz type spaces and related function spaces have interesting applications in studying boundedness of many operators and the theory of partial differential equations; see, for example [5,10,12,14,16].

On the other hand, due to their wide applications in partial dierential equations with non-standard growth ([6]), electrorheological fluids ([13]) and image processing ([1]), the theory of function spaces with variable exponents have attracted a lot of attentions in recent years. Particularly, such theory have achieved great progresses after the notable work of Kov\'{a}\v{c}ik and R\'{a}kosn\'{i}k
[11] in 1991; see [3], [4] and the references therein. In 2010, Izuki [9] studied the Herz spaces with variable exponent $\dot{K}^{\alpha,p}_{q(\cdot)}(\mathbb{R}^{n})$ ($K^{\alpha,p}_{q(\cdot)}(\mathbb{R}^{n})$) and proved the boundedness of some sublinear operators on the spaces. In addition, Wang and Liu [15] introduced a certain Herz-type Hardy spaces with variable exponent $H\dot{K}^{\alpha,p}_{q(\cdot)}(\mathbb{R}^{n})$ ($HK^{\alpha,p}_{q(\cdot)}(\mathbb{R}^{n})$) in 2012. Recently, Yan et al. introduced the variable weak Hardy spaces $WH^{p(\cdot)}(\mathbb{R}^{n})$ and established some characterizations and the boundedness of some operators in [17].

Motivated by [8,9,15,17], we aim to develop the weak Herz spaces and the weak Herz-type Hardy spaces in the setting of variable exponents. In Section 2, we first briefly recall some standard notations and lemmas in variable function spaces. Then we will define the weak Herz spaces and the weak Herz-type Hardy spaces with variable exponent. In Section 3, firstly we will prove the boundedness for a class of sublinear operators from $\dot{K}^{\alpha,p}_{q(\cdot)}(\mathbb{R}^{n})$ (or $K^{\alpha,p}_{q(\cdot)}(\mathbb{R}^{n})$) into the weak Herz spaces with variable exponent $W\dot{K}^{\alpha,p}_{q(\cdot)}(\mathbb{R}^{n})$ (or $WK^{\alpha,p}_{q(\cdot)}(\mathbb{R}^{n})$). Then we will establish the boundedness of local Calder\'{o}n-Zygmund type operators from $H\dot{K}^{\alpha,p}_{q(\cdot)}(\mathbb{R}^{n})$ (or $HK^{\alpha,p}_{q(\cdot)}(\mathbb{R}^{n})$) into $\dot{K}^{\alpha,p}_{q(\cdot)}(\mathbb{R}^{n})$ (or $K^{\alpha,p}_{q(\cdot)}(\mathbb{R}^{n})$) and from $H\dot{K}^{\alpha,p}_{q(\cdot)}(\mathbb{R}^{n})$ (or $HK^{\alpha,p}_{q(\cdot)}(\mathbb{R}^{n})$) into $W\dot{K}^{\alpha,p}_{q(\cdot)}(\mathbb{R}^{n})$ (or $WK^{\alpha,p}_{q(\cdot)}(\mathbb{R}^{n})$) at the critical index. Subsequently, we will show the standard Calder\'{o}n-Zygmund operators is bounded on $H\dot{K}^{\alpha,p}_{q(\cdot)}(\mathbb{R}^{n})$ (or $HK^{\alpha,p}_{q(\cdot)}(\mathbb{R}^{n})$) and from $H\dot{K}^{\alpha,p}_{q(\cdot)}(\mathbb{R}^{n})$ (or $HK^{\alpha,p}_{q(\cdot)}(\mathbb{R}^{n})$) into the weak Herz-type Hardy spaces with variable exponent $WH\dot{K}^{\alpha,p}_{q(\cdot)}(\mathbb{R}^{n})$ (or $WHK^{\alpha,p}_{q(\cdot)}(\mathbb{R}^{n})$) at the endpoint case.

In addition, we denote the Lebesgue measure and the characteristic
function of a measurable set $A\subset \mathbb{R}^{n}$ by $|A|$ and
$\chi_A$ respectively. The notation $f\approx g$ means that there exist constants $C_1,C_2>0$ such that $C_1g\leq f\leq C_2g$.

\section{Weak Herz type spaces with variable exponent\label{s2}}\indent

In this section we first recall some basic definitions and properties of variable function spaces, and then introduce weak Herz spaces and weak Herz-type Hardy spaces with variable exponent.

Given an open set $E\subset \mathbb{R}^{n}$, and a measurable
function $p(\cdot):E\longrightarrow[1,\infty),$
$L^{p(\cdot)}(E)$ denotes the set of measurable functions $f$
on $E$ such that for some $\eta>0,$
$$\int_E\left(\frac{|f(x)|}{\eta}\right)^{p(x)}dx < \infty.$$
This set becomes a Banach function space when equipped with the
Luxemburg-Nakano norm
$$\|f\|_{L^{p(\cdot)}(E)}=\inf\left\{\eta>0:\int_E
\left(\frac{|f(x)|}{\eta}\right)^{p(x)}dx \leq 1\right\}.$$ These
spaces are referred as variable $L^{p}$ spaces, since they
generalized the standard $L^{p}$ spaces: if $p(x)=p$ is constant,
then $L^{p(\cdot)}(E)$ is isometrically isomorphic to
$L^{p}(E)$.

The space $L_{\rm
loc}^{p(\cdot)}(E)$ is defined by $$L_{\rm
loc}^{p(\cdot)}(E):=\{f: f\in L^{p(\cdot)}(F) \,\mathrm{\,for \,\,all \,\,compact \,\,subsets} \,\,F\subset E\}.$$

Define
$\mathcal{P}(E)$ to be set of measurable functions
$p(\cdot):E\longrightarrow[1,\infty)$ such that
$$p^{-}=\mathrm{ess} \inf\{p(x):x\in E\}>1,\quad  p^{+}=\mathrm{ess} \sup\{p(x):x\in E\}<\infty.$$
Denote $p'(x)=p(x)/(p(x)-1).$

Let $f$ be a locally integrable function. The Hardy-Littlewood
maximal operator is defined by
$$Mf(x)=\sup_{B\ni x}\frac{1}{|B|}\int_{B\cap E}|f(y)|dy,$$
where the supremum is taken over all balls $B$ containing $x$. Let $\mathcal{B}(E)$ be the set of $p(\cdot)\in$
$\mathcal{P}(E)$
 such that $M$ is bounded on $L^{p(\cdot)}(E)$.

In variable $L^{p}$ spaces there are some important lemmas as
follows.

\noindent{\bf{Lemma 2.1}}$^{[2]}$\quad \textit{Given an open set $E\subset \mathbb{R}^{n}$. If
$p(\cdot)\in$ $\mathcal{P}(E)$ and satisfies
$$ |p(x)-p(y)|\leq\frac{C}{-\log(|x-y|)}, \quad |x-y|\leq1/2  \eqno(2.1)$$
and
$$ |p(x)-p(y)|\leq\frac{C}{\log(|x|+e)}, \quad |y|\geq|x|, \eqno(2.2)$$
then $p(\cdot)\in$ $\mathcal{B}(E)$, that is the
Hardy-Littlewood maximal operator $M$ is bounded on
$L^{p(\cdot)}(E)$.}

\noindent{\bf{Lemma 2.2}}$^{[11]}$(Generalized H\"{o}lder's inequality)\quad \textit{Given an open set $E\subset \mathbb{R}^{n}$ and let $p(\cdot)\in
\mathcal{P}(E)$. If $f\in L^{p(\cdot)}(E)$ and $g\in
L^{p'(\cdot)}(E)$, then $fg$ is integrable on $E$ and
$$\int_{E}|f(x)g(x)|dx \leq r_{p}\|f\|_{L^{p(\cdot)}(E)}\|g\|_{L^{p'(\cdot)}(E)},$$
where $$r_{p}=1+1/p^--1/p^+.$$}

\noindent{\bf{Lemma 2.3}}$^{[9]}$\quad \textit{Let $p(\cdot)\in
\mathcal{B}(\mathbb{R}^{n})$. Then there exists a positive constant
$C$ such that for all balls $B$ in $\mathbb{R}^{n}$ and all
measurable subsets $S\subset B$,}

$\displaystyle\frac{\|\chi_B\|_{L^{p(\cdot)}(\mathbb{R}^{n})}}{\|\chi_S\|_{L^{p(\cdot)}(\mathbb{R}^{n})}}\leq
C
\frac{|B|}{|S|},\quad\frac{\|\chi_S\|_{L^{p(\cdot)}(\mathbb{R}^{n})}}{\|\chi_B\|_{L^{p(\cdot)}(\mathbb{R}^{n})}}\leq
C \left(\frac{|S|}{|B|}\right)^{\delta_1}$ \textit{and}
$\displaystyle\frac{\|\chi_S\|_{L^{p'(\cdot)}(\mathbb{R}^{n})}}{\|\chi_B\|_{L^{p'(\cdot)}(\mathbb{R}^{n})}}\leq
C \left(\frac{|S|}{|B|}\right)^{\delta_2},$

\noindent \textit{where $0<\delta_1,
\delta_2<1$ depend on $p(\cdot)$.}

Throughout this paper $\delta_2$ is the same as in Lemma 2.3.

\noindent{\bf{Lemma 2.4}}$^{[9]}$\quad \textit{Suppose $p(\cdot)\in
\mathcal{B}(\mathbb{R}^{n})$. Then there exists a positive constant $C$
such that for all balls $B$ in $\mathbb{R}^{n}$,
$$\frac{1}{|B|}\|\chi_B\|_{L^{p(\cdot)}(\mathbb{R}^{n})}\|\chi_B\|_{L^{p'(\cdot)}(\mathbb{R}^{n})}\leq
C.$$}

\noindent{\bf Lemma 2.5}$^{[3]}$\quad \textit{Given $E$ and $p(\cdot)\in \mathcal{P}(E)$, let $f:E\times E\rightarrow\mathbb{R}$ be a measurable
function (with respect to product measure) such that for almost every $y\in E,\,f(\cdot,y)\in L^{p(\cdot)}(E)$. Then
$$\left\|\int_{E}f(\cdot,y)dy\right\|_{L^{p(\cdot)}(E)}\leq C\int_{E}\left\|f(\cdot,y)\right\|_{L^{p(\cdot)}(E)}dy.$$}

Now we recall the definition of the Herz-type spaces with variable
exponent. Let $B_k=\{x\in \mathbb{R}^{n}: |x|\leq 2^k\}$ and
$A_k=B_k\setminus B_{k-1}$ for $k\in \mathbb{Z}$. Denote
$\mathbb{Z_+}$ and $\mathbb{N}$ as the sets of all positive and
non-negative integers, $\chi_k=\chi_{A_k}$ for $k\in \mathbb{Z}$,
$\tilde{\chi}_k=\chi_k$ if $k\in \mathbb{Z_+}$ and
$\tilde{\chi}_0=\chi_{B_0}$.

\noindent {\bf Definition 2.1}$^{[9]}$\quad Let $\alpha\in
\mathbb{R}, 0<p\leq \infty$ and $q(\cdot)\in
\mathcal{P}(\mathbb{R}^{n})$. The homogeneous Herz space with
variable exponent $\dot{K}^{\alpha,p}_{q(\cdot)}(\mathbb{R}^{n})$ is
defined by
$$\dot{K}^{\alpha,p}_{q(\cdot)}(\mathbb{R}^{n})=\{f\in L_{\rm
loc}^{q(\cdot)}(\mathbb{R}^{n}\setminus \{0\}):
\|f\|_{\dot{K}^{\alpha,p}_{q(\cdot)}(\mathbb{R}^{n})}<\infty\},$$
where
$$\|f\|_{\dot{K}^{\alpha,p}_{q(\cdot)}(\mathbb{R}^{n})}=\left\{\sum_{k=-\infty}^\infty2^{k\alpha
p}\|f\chi_k\|^p_{L^{q(\cdot)}(\mathbb{R}^{n})}\right\}^{1/p}.$$ The
non-homogeneous Herz space with variable exponent
$K^{\alpha,p}_{q(\cdot)}(\mathbb{R}^{n})$ is defined by
$$K^{\alpha,p}_{q(\cdot)}(\mathbb{R}^{n})=\{f\in L_{\rm
loc}^{q(\cdot)}(\mathbb{R}^{n}):
\|f\|_{K^{\alpha,p}_{q(\cdot)}(\mathbb{R}^{n})}<\infty\},$$ where
$$\|f\|_{K^{\alpha,p}_{q(\cdot)}(\mathbb{R}^{n})}=\left\{\sum_{k=0}^\infty2^{k\alpha
p}\|f\tilde{\chi}_k\|^p_{L^{q(\cdot)}(\mathbb{R}^{n})}\right\}^{1/p}.$$

In [15], the authors gave the definition of Herz-type Hardy space
with variable exponent
$H\dot{K}^{\alpha,p}_{q(\cdot)}(\mathbb{R}^{n})$ and the atomic
decomposition characterizations. $\mathcal{S}(\mathbb{R}^{n})$
denotes the space of Schwartz functions, and
$\mathcal{S'}(\mathbb{R}^{n})$ denotes the dual space of
$\mathcal{S}(\mathbb{R}^{n})$. Let $G_N(f)(x)$ be the grand maximal
function of $f(x)$ defined by
$$G_N(f)(x)=\sup_{\phi\in\mathcal{A}_N}|\phi^*_\nabla(f)(x)|,$$
where
$\displaystyle\mathcal{A}_N=\{\phi\in\mathcal{S}(\mathbb{R}^{n}):\sup_{|\alpha|,|\beta|\leq
N}|x^\alpha D^\beta\phi(x)|\leq 1\}$ and $N>n+1$, $\phi^*_\nabla$ is
the nontangential maximal operator defined by
$$\phi^*_\nabla(f)(x)=\sup_{|y-x|<t}|\phi_t\ast f(y)|$$
with $\phi_t(x)=t^{-n}\phi(x/t)$.

\noindent{\bf Definition 2.2}$^{[15]}$\quad Let $\alpha\in\mathbb{R},
0<p<\infty, q(\cdot)\in \mathcal{P}(\mathbb{R}^{n})$ and $N>n+1$.

(i) The homogeneous Herz-type Hardy space with variable exponent
$H\dot{K}^{\alpha,p}_{q(\cdot)}(\mathbb{R}^{n})$ is defined by
$$H\dot{K}^{\alpha,p}_{q(\cdot)}(\mathbb{R}^{n})=\left\{f\in \mathcal{S'}(\mathbb{R}^{n}):
G_N(f)(x)\in\dot{K}^{\alpha,p}_{q(\cdot)}(\mathbb{R}^{n})\right\}$$ and
$\|f\|_{H\dot{K}^{\alpha,p}_{q(\cdot)}(\mathbb{R}^{n})}=\|G_N(f)\|_{\dot{K}^{\alpha,p}_{q(\cdot)}(\mathbb{R}^{n})}$.

(ii) The non-homogeneous Herz-type Hardy space with variable exponent
$HK^{\alpha,p}_{q(\cdot)}(\mathbb{R}^{n})$ is defined by
$$HK^{\alpha,p}_{q(\cdot)}(\mathbb{R}^{n})=\left\{f\in \mathcal{S'}(\mathbb{R}^{n}):
G_N(f)(x)\in K^{\alpha,p}_{q(\cdot)}(\mathbb{R}^{n})\right\}$$ and
$\|f\|_{HK^{\alpha,p}_{q(\cdot)}(\mathbb{R}^{n})}=\|G_N(f)\|_{K^{\alpha,p}_{q(\cdot)}(\mathbb{R}^{n})}$.

For $x\in\mathbb{R}$ we denote by $[x]$ the largest integer less
than or equal to $x$.

\noindent{\bf Definition 2.3}$^{[15]}$\quad Let
$n\delta_2\leq\alpha<\infty, q(\cdot)\in
\mathcal{P}(\mathbb{R}^{n})$, and non-negative integer $s\geq
[\alpha-n\delta_2]$.

(i) A function $a(x)$ on $\mathbb{R}^{n}$ is said to be a central
$(\alpha, q(\cdot))$-atom, if it satisfies

\hspace{3mm}(1) supp\,$a\subset B(0,r)=\{x\in
\mathbb{R}^{n}:|x|<r\}$.

\hspace{3mm}(2) $\|a\|_{L^{q(\cdot)}(\mathbb{R}^{n})}\leq
|B(0,r)|^{-\alpha/n}$.

\hspace{3mm}(3) $\int_{\mathbb{R}^{n}}a(x)x^\beta dx=0, |\beta|\leq
s$.

(ii) A function $a(x)$ on $\mathbb{R}^{n}$ is said to be a central
$(\alpha, q(\cdot))$-atom of restricted type, if it satisfies the
conditions (2), (3) above and

\hspace{3mm}(1)$'$ supp\,$a\subset B(0,r), r\geq 1$.

If $r=2^k$ for some $k\in\mathbb{Z}$ in Definition 1.3, then the
corresponding central $(\alpha, q(\cdot))$-atom is called a dyadic central $(\alpha, q(\cdot))$-atom.

\noindent{\bf Lemma 2.6}$^{[15]}$\quad \textit{Let
$n\delta_2\leq\alpha<\infty, \,\,0<p<\infty$ and $q(\cdot)\in
\mathcal{B}(\mathbb{R}^{n})$. Then $f\in
H\dot{K}^{\alpha,p}_{q(\cdot)}(\mathbb{R}^{n})$(or
$HK^{\alpha,p}_{q(\cdot)}(\mathbb{R}^{n})$) if and only if
$$f=\sum_{k=-\infty}^\infty\lambda_ka_k\,\left(\mathrm{or}\,\sum_{k=0}^\infty\lambda_ka_k\right),\quad \mathrm{in\,\, the\,\, sense\,\, of\,\,} \mathcal{S'}(\mathbb{R}^{n}),$$ where each
$a_k$ is a central $(\alpha, q(\cdot))$-atom(or central $(\alpha,
q(\cdot))$-atom of restricted type) with support contained in $B_k$
and $\displaystyle\sum_{k=-\infty}^\infty|\lambda_k|^p<\infty$(or
$\displaystyle\sum_{k=0}^\infty|\lambda_k|^p<\infty$). Moreover,
$$\|f\|_{H\dot{K}^{\alpha,p}_{q(\cdot)}(\mathbb{R}^{n})}\approx
\inf\left(\sum_{k=-\infty}^\infty|\lambda_k|^p\right)^{1/p}\,\left(\mathrm{or}\quad\|f\|_{HK^{\alpha,p}_{q(\cdot)}(\mathbb{R}^{n})}\approx\inf\left(\sum_{k=0}^\infty|\lambda_k|^p\right)^{1/p}\right),$$
where the infimum is taken over all above decompositions of $f$.}

Next we recall the definition of the weak Herz spaces and the weak Herz-type Hardy spaces. For $k\in \mathbb{Z}$, let $m_k(\sigma, f)=|\{x\in A_k: |f(x)|>\sigma\}|$; for $k\in \mathbb{N}$, let $\tilde{m}_k(\sigma, f)=m_k(\sigma, f)$ and $\tilde{m}_0(\sigma, f)=|\{x\in B(0, 1): |f(x)|>\sigma\}|$.

\noindent {\bf Definition 2.4}$^{[7]}$\quad Let $\alpha\in
\mathbb{R}, 0<p\leq \infty$ and $0<q<\infty$. A measurable function $f(x)$ on $\mathbb{R}^{n}$ is said to belong to the homogeneous weak Herz space
$W\dot{K}^{\alpha,p}_{q}(\mathbb{R}^{n})$, if
$$\|f\|_{W\dot{K}^{\alpha,p}_{q}(\mathbb{R}^{n})}=\sup_{\lambda>0}\lambda\left\{\sum_{k=-\infty}^\infty2^{k\alpha
p}m_k(\lambda, f)^{p/q}\right\}^{1/p}<\infty,$$
where the usual modification is made when $p=\infty$.

A measurable function $f(x)$ on $\mathbb{R}^{n}$ is said to belong to the non-homogeneous weak Herz space
$WK^{\alpha,p}_{q}(\mathbb{R}^{n})$, if
$$\|f\|_{WK^{\alpha,p}_{q}(\mathbb{R}^{n})}=\sup_{\lambda>0}\lambda\left\{\sum_{k=0}^\infty2^{k\alpha
p}\tilde{m}_k(\lambda, f)^{p/q}\right\}^{1/p}<\infty,$$
where the usual modification is made when $p=\infty$.

\noindent{\bf Definition 2.5}$^{[8]}$\quad Let $\alpha\in\mathbb{R},
0<p, q\leq\infty$ and $N$ be sufficiently large. We define the spaces
$$WH\dot{K}^{\alpha,p}_{q}(\mathbb{R}^{n})=\{f\in \mathcal{S'}(\mathbb{R}^{n}):
G_N(f)(x)\in W\dot{K}^{\alpha,p}_{q}(\mathbb{R}^{n})\}$$ and
$$WHK^{\alpha,p}_{q}(\mathbb{R}^{n})=\{f\in \mathcal{S'}(\mathbb{R}^{n}):
G_N(f)(x)\in WK^{\alpha,p}_{q}(\mathbb{R}^{n})\}.$$
Moreover, we define that

$\|f\|_{WH\dot{K}^{\alpha,p}_{q}(\mathbb{R}^{n})}=\|G_N(f)\|_{W\dot{K}^{\alpha,p}_{q}(\mathbb{R}^{n})}$ and
$\|f\|_{WHK^{\alpha,p}_{q}(\mathbb{R}^{n})}=\|G_N(f)\|_{WK^{\alpha,p}_{q}(\mathbb{R}^{n})}$.

Now we extend the Definitions 2.4 and 2.5 to the case of function spaces with variable exponent.

\noindent {\bf Definition 2.6}\quad Let $\alpha\in
\mathbb{R}, 0<p\leq \infty$ and $q(\cdot)\in
\mathcal{P}(\mathbb{R}^{n})$. A measurable function $f(x)$ on $\mathbb{R}^{n}$ is said to belong to the homogeneous weak Herz space with variable exponent
$W\dot{K}^{\alpha,p}_{q(\cdot)}(\mathbb{R}^{n})$, if
$$\|f\|_{W\dot{K}^{\alpha,p}_{q(\cdot)}(\mathbb{R}^{n})}=\sup_{\lambda>0}\lambda\left\{\sum_{k=-\infty}^\infty2^{k\alpha
p}\|\chi_{\{x\in A_k: |f(x)|>\lambda\}}\|_{L^{q(\cdot)}(\mathbb{R}^{n})}^p\right\}^{1/p}<\infty,$$
where the usual modification is made when $p=\infty$.

A measurable function $f(x)$ on $\mathbb{R}^{n}$ is said to belong to the non-homogeneous weak Herz space with variable exponent
$WK^{\alpha,p}_{q(\cdot)}(\mathbb{R}^{n})$, if
$$\|f\|_{WK^{\alpha,p}_{q(\cdot)}(\mathbb{R}^{n})}=\sup_{\lambda>0}\lambda\left\{\sum_{k=0}^\infty2^{k\alpha
p}\|\chi_{\{x\in A_k: |f(x)|>\lambda\}}\|_{L^{q(\cdot)}(\mathbb{R}^{n})}^p\right\}^{1/p}<\infty,$$
where the usual modification is made when $p=\infty$.

\noindent {\bf Definition 2.7}\quad Let $\alpha\in\mathbb{R},
0<p<\infty, q(\cdot)\in \mathcal{P}(\mathbb{R}^{n})$ and $N>n+1$. We define the spaces
$$WH\dot{K}^{\alpha,p}_{q(\cdot)}(\mathbb{R}^{n})=\{f\in \mathcal{S'}(\mathbb{R}^{n}):
G_N(f)(x)\in W\dot{K}^{\alpha,p}_{q(\cdot)}(\mathbb{R}^{n})\}$$ and
$$WHK^{\alpha,p}_{q(\cdot)}(\mathbb{R}^{n})=\{f\in \mathcal{S'}(\mathbb{R}^{n}):
G_N(f)(x)\in WK^{\alpha,p}_{q(\cdot)}(\mathbb{R}^{n})\}.$$
Moreover, we define that

$\|f\|_{WH\dot{K}^{\alpha,p}_{q(\cdot)}(\mathbb{R}^{n})}=\|G_N(f)\|_{W\dot{K}^{\alpha,p}_{q(\cdot)}(\mathbb{R}^{n})}$ and
$\|f\|_{WHK^{\alpha,p}_{q(\cdot)}(\mathbb{R}^{n})}=\|G_N(f)\|_{WK^{\alpha,p}_{q(\cdot)}(\mathbb{R}^{n})}$.

\noindent {\bf Remark 2.1}\quad If $q(\cdot)=q$ is constant in Definition 2.6 and Definition 2.7, then we can easily get the Definition 2.4 and Definition 2.5 respectively.

\noindent {\bf Remark 2.2}\quad If $\alpha=0$, then $WL^{q(\cdot)}(\mathbb{R}^{n})\subset W\dot{K}^{0,p}_{q(\cdot)}(\mathbb{R}^{n})$, where $WL^{q(\cdot)}(\mathbb{R}^{n})$ is the weak Lebesgue space with variable exponent and $$\|f\|_{WL^{q(\cdot)}(\mathbb{R}^{n})}=\sup_{\lambda>0}\lambda
\|\chi_{\{x\in A_k: |f(x)|>\lambda\}}\|_{L^{q(\cdot)}(\mathbb{R}^{n})}<\infty.$$

\section{Boundedness of some singular integral operators\label{s3}}\indent

Firstly we establish the boundedness for a large class of operators on Herz type spaces with variable exponent.

\noindent{\bf{Theorem 3.1}}\quad \textit{Let $0<\alpha<n\delta_2$, $q(\cdot)\in
\mathcal{B}(\mathbb{R}^{n})$, $p\in (0,1]$. If $T$ be a sublinear operator and bounded on $WL^{q(\cdot)}(\mathbb{R}^{n})$ satisfying
$$|Tf(x)|\leq C|x|^{-n}\int_{\mathbb{R}^{n}}|f(y)|dy \eqno(3.1)$$
for any $f\in L_{\rm
loc}^{1}(\mathbb{R}^{n})$, supp\,$f\subset B(0,r)$ and $x\notin B(0,2r)$.
Then $T$ maps continuously $\dot{K}^{\alpha,p}_{q(\cdot)}(\mathbb{R}^{n})$ (or $K^{\alpha,p}_{q(\cdot)}(\mathbb{R}^{n})$) into $W\dot{K}^{\alpha,p}_{q(\cdot)}(\mathbb{R}^{n})$ (or $WK^{\alpha,p}_{q(\cdot)}(\mathbb{R}^{n})$).}

\noindent{\bf{Proof}}\quad We only prove the homogeneous case. The
non-homogeneous case can be proved in the same way. For any $k\in\mathbb{Z}$, we decompose $f$ into
$$f(x)=f(x)\chi_{\{|x|\leq2^{k-3}\}}(x)+f(x)\chi_{\{|x|>2^{k-3}\}}(x)=f_1(x)+f_2(x).$$
Then $|Tf(x)|\leq|Tf_1(x)|+|Tf_2(x)|$,
and
$$\begin{array}{rl}
\displaystyle\|Tf\|_{W\dot{K}^{\alpha,p}_{q(\cdot)}(\mathbb{R}^{n})}
&\displaystyle=\sup_{\lambda>0}\lambda\left\{\sum_{k=-\infty}^{\infty}2^{k\alpha p}\|\chi_{\{x\in A_k: |Tf(x)|>\lambda\}}\|^p_{L^{q(\cdot)}(\mathbb{R}^{n})}\right\}^{1/p}\\
&\displaystyle\leq C\sup_{\lambda>0}\lambda\left\{\sum_{k=-\infty}^{\infty}2^{k\alpha p}\|\chi_{\{x\in A_k: |Tf_1(x)|>\frac{\lambda}{2}\}}\|^p_{L^{q(\cdot)}(\mathbb{R}^{n})}\right\}^{1/p}\\

\end{array}$$
$$\begin{array}{rl}
&\displaystyle\hspace{3mm}+C\sup_{\lambda>0}\lambda\left\{\sum_{k=-\infty}^{\infty}2^{k\alpha p}\|\chi_{\{x\in A_k: |Tf_2(x)|>\frac{\lambda}{2}\}}\|^p_{L^{q(\cdot)}(\mathbb{R}^{n})}\right\}^{1/p}\\
&\displaystyle=I_1+I_2.

\end{array}$$
Using the $WL^{q(\cdot)}(\mathbb{R}^{n})$-boundedness of $T$ and $0<p\leq 1$, we have
$$\begin{array}{rl}
\displaystyle I_2&\displaystyle
\leq C\left\{\sum_{k=-\infty}^{\infty}2^{k\alpha p}\|Tf_2\chi_k\|^p_{WL^{q(\cdot)}(\mathbb{R}^{n})}\right\}^{1/p}\\
&\displaystyle\leq C\left\{\sum_{k=-\infty}^{\infty}2^{k\alpha p}\left(\sum_{l=k-2}^{\infty}\|f\chi_l\|_{L^{q(\cdot)}(\mathbb{R}^{n})}\right)^p\right\}^{1/p}\\
&\displaystyle\leq C\left\{\sum_{l=-\infty}^{\infty}\|f\chi_l\|^p_{L^{q(\cdot)}(\mathbb{R}^{n})}\left(\sum_{k=-\infty}^{l+2}2^{k\alpha p}\right)\right\}^{1/p}\\
&\displaystyle\leq C\left\{\sum_{l=-\infty}^{\infty}2^{l\alpha p}\|f\chi_l\|^p_{L^{q(\cdot)}(\mathbb{R}^{n})}\right\}^{1/p}\\
&\displaystyle=C\|f\|_{\dot{K}^{\alpha,p}_{q(\cdot)}(\mathbb{R}^{n})}.

\end{array}$$

For $I_1$, noting that $x\in A_k$ and supp\,$f_1\subset\{x\in \mathbb{R}^{n}: |x|\leq 2^{k-3}\}$, by (3.1) and the generalized
H\"{o}lder inequality we have
$$\begin{array}{rl}
\displaystyle |Tf_1(x)|&\displaystyle
\leq C|x|^{-n}\|f_1\|_{L^1(\mathbb{R}^{n})}\\
&\displaystyle\leq C2^{-kn}\sum_{l=-\infty}^{k-3}\|f\chi_l\|_{L^{1}(\mathbb{R}^{n})}\\
&\displaystyle\leq C2^{-kn}\sum_{l=-\infty}^{k-3}\|f\chi_l\|_{L^{q(\cdot)}(\mathbb{R}^{n})}\|\chi_l\|_{L^{q'(\cdot)}(\mathbb{R}^{n})}.

\end{array}$$
So by Lemmas 2.3 and 2.4 we have
$$\begin{array}{rl}
&\displaystyle \left\|\chi_{\{x\in A_k: |Tf_1(x)|>\frac{\lambda}{2}\}}\right\|_{L^{q(\cdot)}(\mathbb{R}^{n})}\\
&\displaystyle\leq \left\|\chi_{\{x\in A_k: C2^{-kn}\sum_{l=-\infty}^{k-3}\|f\chi_l\|_{L^{q(\cdot)}(\mathbb{R}^{n})}\|\chi_l\|_{L^{q'(\cdot)}(\mathbb{R}^{n})}>\frac{\lambda}{2}\}}\right\|_{L^{q(\cdot)}(\mathbb{R}^{n})}\\
&\displaystyle\leq C2^{-kn+1}\lambda^{-1}\sum_{l=-\infty}^{k-3}\|f\chi_l\|_{L^{q(\cdot)}(\mathbb{R}^{n})}\|\chi_l\|_{L^{q'(\cdot)}(\mathbb{R}^{n})}\|\chi_k\|_{L^{q(\cdot)}(\mathbb{R}^{n})}\\
&\displaystyle\leq C2^{-kn}\lambda^{-1}\sum_{l=-\infty}^{k-3}\|f\chi_l\|_{L^{q(\cdot)}(\mathbb{R}^{n})}\|\chi_{B_l}\|_{L^{q'(\cdot)}(\mathbb{R}^{n})}\|\chi_{B_k}\|^{-1}_{L^{q'(\cdot)}(\mathbb{R}^{n})}|B_k|\\
&\displaystyle\leq C\lambda^{-1}\sum_{l=-\infty}^{k-3}\|f\chi_l\|_{L^{q(\cdot)}(\mathbb{R}^{n})}2^{(l-k)n\delta_2}.

\end{array}$$
Thus by $0<p\leq1$ and $0<\alpha<n\delta_2$ we obtain
$$\begin{array}{rl}
\displaystyle I_1&\displaystyle
=C\sup_{\lambda>0}\lambda\left\{\sum_{k=-\infty}^{\infty}2^{k\alpha p}\|\chi_{\{x\in A_k: |Tf_1(x)|>\frac{\lambda}{2}\}}\|^p_{L^{q(\cdot)}(\mathbb{R}^{n})}\right\}^{1/p}\\
&\displaystyle\leq C\left\{\sum_{k=-\infty}^{\infty}2^{k\alpha p}\sum_{l=-\infty}^{k-3}\|f\chi_l\|^p_{L^{q(\cdot)}(\mathbb{R}^{n})}2^{(l-k)n\delta_2p}\right\}^{1/p}\\

\end{array}$$
$$\begin{array}{rl}
&\displaystyle\leq C\left\{\sum_{l=-\infty}^{\infty}2^{l\alpha p}\|f\chi_l\|^p_{L^{q(\cdot)}(\mathbb{R}^{n})}\left(\sum_{k=l+3}^{\infty}2^{(l-k)(n\delta_2-\alpha) p}\right)\right\}^{1/p}\\
&\displaystyle\leq C\left\{\sum_{l=-\infty}^{\infty}2^{l\alpha p}\|f\chi_l\|^p_{L^{q(\cdot)}(\mathbb{R}^{n})}\right\}^{1/p}\\
&\displaystyle=C\|f\|_{\dot{K}^{\alpha,p}_{q(\cdot)}(\mathbb{R}^{n})}.

\end{array}$$

This completes the proof of Theorem 3.1.

\noindent {\bf Remark 3.1}\quad If the condition (3.1) is substituted by the following condition
$$|Tf(x)|\leq C\int_{\mathbb{R}^{n}}\frac{|f(y)|}{|x-y|^{n}}dy, \quad x\notin\mathrm{supp}f, \eqno(3.2)$$
then the conclusions of Theorem 3.1 are still true.

\noindent {\bf Remark 3.2}\quad The condition (3.1) is very weak and many classical operators satisfy (3.1), such as Calder\'{o}n-Zygmund operators, multipliers and oscillatory singular integrals, Bochner-Riesz operators at the critical index and so on.

Now we turn our attention to the behaviours of local Calder\'{o}n-Zygmund type operators on the Herz-type Hardy spaces with variable exponent.

\noindent{\bf{Theorem 3.2}}\quad \textit{Let $T: \mathcal{S}(\mathbb{R}^{n})\rightarrow\mathcal{S'}(\mathbb{R}^{n})$ be a linear and continuous operator. Assume that the distribution kernel of $T$ coincides in the complement of the diagonal with a locally integrable function $k(x,y)$ which satisfies
$$\sup_{y\in B_k}\|[k(\cdot, y)-k(\cdot, 0)]\chi_l\|_{L^{q(\cdot)}(\mathbb{R}^{n})}\leq C2^{(k-l)\delta-ln}\|\chi_l\|_{L^{q(\cdot)}(\mathbb{R}^{n})}\eqno(3.3)$$
for $k, l\in\mathbb{Z}$, some $\delta\in(0,1]$ and some $q(\cdot)\in
\mathcal{B}(\mathbb{R}^{n})$. Suppose that $T$ can be extended to a continuous operator on $L^{q(\cdot)}(\mathbb{R}^{n})$. If $n\delta_2\leq\alpha<n\delta_2+\delta$ and $0<p\leq \infty$, then $T$ maps continuously $H\dot{K}^{\alpha,p}_{q(\cdot)}(\mathbb{R}^{n})$ into
$\dot{K}^{\alpha,p}_{q(\cdot)}(\mathbb{R}^{n})$.}

\noindent {\bf Remark 3.3}\quad If the condition (3.3) is true only for $k\in\mathbb{N}$ and $l\in \mathbb{Z_+}$, then the operator $T$ in Theorem 3.2 maps continously $HK^{\alpha,p}_{q(\cdot)}(\mathbb{R}^{n})$ into $K^{\alpha,p}_{q(\cdot)}(\mathbb{R}^{n})$.

\noindent {\bf Remark 3.4}\quad Assuming more regularity on the kernel $k(x,y)$, we can extend Theorem 3.2 to larger range of $\alpha$.

\noindent{\bf Proof of Theorem 3.2}\quad Firstly we suppose that $0<p\leq 1$. In this case, we only need to prove $\|Ta_k\|_{\dot{K}^{\alpha,p}_{q(\cdot)}(\mathbb{R}^{n})}\leq C$, where $a_k$ is a dyadic
central $(\alpha, q(\cdot))$-atom with the support $B_k$ and $C$ is independent of $k$.
Write
$$\begin{array}{rl}
\displaystyle\|Ta_k\|_{\dot{K}^{\alpha,p}_{q(\cdot)}(\mathbb{R}^{n})}&\displaystyle\leq C\left\{\sum_{l=-\infty}^{k+3}2^{l\alpha p}\|(Ta_k)\chi_l\|^p_{L^{q(\cdot)}(\mathbb{R}^{n})}\right\}^{1/p}\\
&\displaystyle\hspace{3mm}+C\left\{\sum_{l=k+4}^{\infty}2^{l\alpha p}\|(Ta_k)\chi_l\|^p_{L^{q(\cdot)}(\mathbb{R}^{n})}\right\}^{1/p}\\
&\displaystyle=J_1+J_2.

\end{array}$$
Since $T$ is bounded on $L^{q(\cdot)}(\mathbb{R}^{n})$, by Definition 2.3 we have
$$\begin{array}{rl}
\displaystyle J_1&\displaystyle
\leq C\left\{\sum_{l=-\infty}^{k+3}2^{l\alpha p}\|Ta_k\|^p_{L^{q(\cdot)}(\mathbb{R}^{n})}\right\}^{1/p}\\

\end{array}$$
$$\begin{array}{rl}
&\displaystyle\leq C\left\{\sum_{l=-\infty}^{k+3}2^{l\alpha p}\|a_k\|^p_{L^{q(\cdot)}(\mathbb{R}^{n})}\right\}^{1/p}\\
&\displaystyle\leq C 2^{-k\alpha}\left\{\sum_{l=-\infty}^{k+3}2^{l\alpha p}\right\}^{1/p}\leq C,

\end{array}$$
where $C$ is independent of $k$.

For $J_2$, using (3.3), Lemmas 2.3, 2.4, 2.5 and the generalized
H\"{o}lder inequality we have
$$\begin{array}{rl}
\displaystyle \|(Ta_k)\chi_l\|_{L^{q(\cdot)}(\mathbb{R}^{n})}&\displaystyle
\leq C\int_{B_k}\|[k(\cdot, y)-k(\cdot, 0)]\chi_l\|_{L^{q(\cdot)}(\mathbb{R}^{n})}|a_k(y)|dy\\
&\displaystyle\leq C2^{(k-l)\delta-ln}\|\chi_l\|_{L^{q(\cdot)}(\mathbb{R}^{n})}\|a_k\|_{L^{q(\cdot)}(\mathbb{R}^{n})}\|\chi_k\|_{L^{q'(\cdot)}(\mathbb{R}^{n})}\\
&\displaystyle\leq C2^{(k-l)\delta-k\alpha}\|\chi_{B_l}\|^{-1}_{L^{q'(\cdot)}(\mathbb{R}^{n})}\|\chi_{B_k}\|_{L^{q'(\cdot)}(\mathbb{R}^{n})}\\
&\displaystyle\leq C2^{(k-l)(\delta+n\delta_2)-k\alpha}.

\end{array}$$
Thus by $\alpha<n\delta_2+\delta$ we obtain
$$J_2\leq C\left\{\sum_{l=k+4}^{\infty}2^{(l-k)(\alpha-\delta-n\delta_2)p}\right\}^{1/p}\leq C,$$
where $C$ is independent of $k$.
This finishes the proof for the case $0<p\leq1$. Now let $1<p<\infty$ and $f\in
H\dot{K}^{\alpha,p}_{q(\cdot)}(\mathbb{R}^{n})$. By Lemma 2.6 we get
$\displaystyle f=\sum_{k=-\infty}^\infty\lambda_ka_k$, where
$\displaystyle\|f\|_{H\dot{K}^{\alpha,p}_{q(\cdot)}(\mathbb{R}^{n})}\approx
\inf(\sum_{k=-\infty}^\infty|\lambda_k|^{p})^{1/p}$ (the infimum
is taken over above decompositions of $f$), and $a_k$ is a dyadic
central $(\alpha, q(\cdot))$-atom with the support $B_k$. Write
$$\begin{array}{rl}
\displaystyle\|Tf\|_{\dot{K}^{\alpha,p}_{q(\cdot)}(\mathbb{R}^{n})}&\displaystyle\leq C\left\{\sum_{l=-\infty}^{\infty}2^{l\alpha p}\left(\sum_{k=-\infty}^{l-4}|\lambda_k|\|(Ta_k)\chi_l\|_{L^{q(\cdot)}(\mathbb{R}^{n})}\right)^p\right\}^{1/p}\\
&\displaystyle\hspace{3mm}+C\left\{\sum_{l=-\infty}^{\infty}2^{l\alpha p}\left(\sum_{k=l-3}^{\infty}|\lambda_k|\|(Ta_k)\chi_l\|_{L^{q(\cdot)}(\mathbb{R}^{n})}\right)^p\right\}^{1/p}\\
&\displaystyle=U_1+U_2.

\end{array}$$
For $U_2$, by the H\"{o}lder inequality and the fact that $T$ is bounded on $L^{q(\cdot)}(\mathbb{R}^{n})$, we have
$$\begin{array}{rl}
\displaystyle U_2&\displaystyle
\leq C\left\{\sum_{l=-\infty}^{\infty}2^{l\alpha p}\left(\sum_{k=l-3}^{\infty}|\lambda_k|\|(Ta_k)\|_{L^{q(\cdot)}(\mathbb{R}^{n})}\right)^p\right\}^{1/p}\\
&\displaystyle\leq C\left\{\sum_{l=-\infty}^{\infty}2^{l\alpha p}\left(\sum_{k=l-3}^{\infty}|\lambda_k|2^{-k\alpha}\right)^p\right\}^{1/p}\\
&\displaystyle\leq C\left\{\sum_{l=-\infty}^{\infty}2^{l\alpha p/2}\left(\sum_{k=l-3}^{\infty}|\lambda_k|^p 2^{-k\alpha p/2}\right)\right\}^{1/p}\\
&\displaystyle\leq C\left\{\sum_{k=-\infty}^{\infty}|\lambda_k|^p\right\}^{1/p}\leq C\|f\|_{H\dot{K}^{\alpha,p}_{q(\cdot)}(\mathbb{R}^{n})}.

\end{array}$$
On the other hand, noting that $\alpha<n\delta_2+\delta$. So by (3.3), Lemmas 2.3, 2.4, 2.5 and the generalized
H\"{o}lder inequality we have
$$\begin{array}{rl}
\displaystyle U_1&\displaystyle
\leq C\left\{\sum_{l=-\infty}^{\infty}2^{l\alpha p}\left(\sum_{k=-\infty}^{l-4}|\lambda_k|\|a_k\|_{L^{1}(\mathbb{R}^{n})}2^{(k-l)\delta-ln}\|\chi_l\|_{L^{q(\cdot)}(\mathbb{R}^{n})}\right)^p\right\}^{1/p}\\
&\displaystyle\leq C\left\{\sum_{l=-\infty}^{\infty}\left(\sum_{k=-\infty}^{l-4}|\lambda_k|2^{(l-k)(\alpha-\delta-n\delta_2)}\right)^p\right\}^{1/p}\\
&\displaystyle\leq C\left\{\sum_{l=-\infty}^{\infty}\left(\sum_{k=-\infty}^{l-4}|\lambda_k|^p2^{(l-k)(\alpha-\delta-n\delta_2)p/2}\right)\right\}^{1/p}\\
&\displaystyle\leq C\left\{\sum_{k=-\infty}^{\infty}|\lambda_k|^p\left(\sum_{l=k+4}^{\infty}2^{(l-k)(\alpha-\delta-n\delta_2)p/2}\right)\right\}^{1/p}\\
&\displaystyle\leq C\left\{\sum_{k=-\infty}^{\infty}|\lambda_k|^p\right\}^{1/p}\leq C\|f\|_{H\dot{K}^{\alpha,p}_{q(\cdot)}(\mathbb{R}^{n})}.

\end{array}$$

Therefore, we complete the proof of Theorem 3.2.

When $\alpha=n\delta_2+\delta$, we have the following result.

\noindent{\bf{Theorem 3.3}}\quad \textit{Let $T: \mathcal{S}(\mathbb{R}^{n})\rightarrow\mathcal{S'}(\mathbb{R}^{n})$ be a linear operator. Suppose that the distribution kernel of $T$ coincides in the complement of the diagonal with a locally integrable function $k(x,y)$ satisfying
$$|k(x, y)-k(x, 0)|\leq C\frac{|y|^\delta}{|x|^{n+\delta}},\eqno(3.4)$$
if $2|y|<|x|$, for some $\delta\in(0,1]$. Suppose that $T$ is bounded on $L^{q(\cdot)}(\mathbb{R}^{n})$ for some $q(\cdot)\in
\mathcal{B}(\mathbb{R}^{n})$, $\alpha=n\delta_2+\delta$ and $0<p\leq 1$, then $T$ maps continuously $H\dot{K}^{\alpha,p}_{q(\cdot)}(\mathbb{R}^{n})$(or
$HK^{\alpha,p}_{q(\cdot)}(\mathbb{R}^{n})$) into
$W\dot{K}^{\alpha,p}_{q(\cdot)}(\mathbb{R}^{n})$(or
$WK^{\alpha,p}_{q(\cdot)}(\mathbb{R}^{n})$).}

\noindent{\bf Proof}\quad We only prove the homogeneous case. Let $\alpha=n\delta_2+\delta$ and $f\in
H\dot{K}^{\alpha,p}_{q(\cdot)}(\mathbb{R}^{n})$. By Lemma 2.6 we get
$\displaystyle f=\sum_{l=-\infty}^\infty\lambda_la_l$, where
$\displaystyle\|f\|_{H\dot{K}^{\alpha,p}_{q(\cdot)}(\mathbb{R}^{n})}\approx
\inf(\sum_{l=-\infty}^\infty|\lambda_l|^{p})^{1/p}$ (the infimum
is taken over above decompositions of $f$), and $a_l$ is a dyadic
central $(\alpha, q(\cdot))$-atom with the support $B_l$. Given $\lambda>0$, we can
write
$$\begin{array}{rl}
&\displaystyle\lambda\left\{\sum_{k=-\infty}^{\infty}2^{k\alpha p}\|\chi_{\{x\in A_k: |Tf(x)|>\lambda\}}\|^p_{L^{q(\cdot)}(\mathbb{R}^{n})}\right\}^{1/p}\\
&\displaystyle\leq C\lambda\left\{\sum_{k=-\infty}^{\infty}2^{k\alpha p}\|\chi_{\{x\in A_k: |\sum_{l=-\infty}^{k-4}\lambda_lTa_l(x)|>\frac{\lambda}{2}\}}\|^p_{L^{q(\cdot)}(\mathbb{R}^{n})}\right\}^{1/p}\\
&\displaystyle\hspace{3mm}+C\lambda\left\{\sum_{k=-\infty}^{\infty}2^{k\alpha p}\|\chi_{\{x\in A_k: |\sum_{l=k-3}^{\infty}\lambda_lTa_l(x)|>\frac{\lambda}{2}\}}\|^p_{L^{q(\cdot)}(\mathbb{R}^{n})}\right\}^{1/p}\\
&\displaystyle=CV_1+CV_2.

\end{array}$$
For $V_2$, by the $L^{q(\cdot)}(\mathbb{R}^{n})$-boundedness of $T$, we have
$$\begin{array}{rl}
\displaystyle V_2&\displaystyle
\leq C\left\{\sum_{k=-\infty}^{\infty}2^{k\alpha p}\left\|\sum_{l=k-3}^{\infty}|\lambda_lTa_l|\chi_k\right\|^p_{L^{q(\cdot)}(\mathbb{R}^{n})}\right\}^{1/p}\\
&\displaystyle\leq C\left\{\sum_{k=-\infty}^{\infty}2^{k\alpha p}\left(\sum_{l=k-3}^{\infty}|\lambda_l|\|a_l\|_{L^{q(\cdot)}(\mathbb{R}^{n})}\right)^p\right\}^{1/p}\\
&\displaystyle\leq C\left\{\sum_{k=-\infty}^{\infty}\sum_{l=k-3}^{\infty}|\lambda_l|^p2^{(k-l)\alpha p}\right\}^{1/p}\\
&\displaystyle\leq C\left\{\sum_{l=-\infty}^{\infty}|\lambda_l|^p\sum_{k=-\infty}^{l+3}2^{(k-l)\alpha p}\right\}^{1/p}\\
&\displaystyle\leq C\left\{\sum_{l=-\infty}^{\infty}|\lambda_l|^p\right\}^{1/p}\leq C\|f\|_{H\dot{K}^{\alpha,p}_{q(\cdot)}(\mathbb{R}^{n})}.

\end{array}$$
To estimate $V_1$, we notice that if $x\in A_k$ and $k\geq l+4$, then by (3.4), $\alpha=n\delta_2+\delta$, Lemma 2.3 and the generalized
H\"{o}lder inequality we have
$$\begin{array}{rl}
\displaystyle |Ta_l(x)|&\displaystyle\leq\int_{\mathbb{R}^{n}}|k(x, y)-k(x, 0)||a_l(y)|dy\\
&\displaystyle\leq C2^{l\delta-k(n+\delta)}\int_{\mathbb{R}^{n}}|a_l(y)|dy\\
&\displaystyle\leq C2^{l\delta-k(n+\delta)}\|a_l\|_{L^{q(\cdot)}(\mathbb{R}^{n})}\|\chi_l\|_{L^{q'(\cdot)}(\mathbb{R}^{n})}\\
&\displaystyle\leq C2^{l(\delta-\alpha)-k(n+\delta)}\|\chi_l\|_{L^{q'(\cdot)}(\mathbb{R}^{n})}\\
&\displaystyle\leq C2^{l(\delta-\alpha)-k(n+\delta)}2^{(l-k)n\delta_2}\|\chi_k\|_{L^{q'(\cdot)}(\mathbb{R}^{n})}\\
&\displaystyle= C2^{-k(n+\delta+n\delta_2)}\|\chi_k\|_{L^{q'(\cdot)}(\mathbb{R}^{n})}.

\end{array}$$
So by $0<p\leq1$ we have
$$\begin{array}{rl}
\displaystyle \left|\sum_{l=-\infty}^{k-4}\lambda_lTa_l(x)\right|
&\displaystyle\leq C2^{-k(n+\delta+n\delta_2)}\|\chi_k\|_{L^{q'(\cdot)}(\mathbb{R}^{n})}\sum_{l=-\infty}^{k-4}|\lambda_l|\\
&\displaystyle\leq C2^{-k(n+\delta+n\delta_2)}\|\chi_k\|_{L^{q'(\cdot)}(\mathbb{R}^{n})}\left(\sum_{l=-\infty}^{k-4}|\lambda_l|^p\right)^{1/p}\\
&\displaystyle\leq C_02^{-k(n+\delta+n\delta_2)}\|\chi_k\|_{L^{q'(\cdot)}(\mathbb{R}^{n})}\|f\|_{H\dot{K}^{\alpha,p}_{q(\cdot)}(\mathbb{R}^{n})}.

\end{array}$$
If $|\{x\in A_k: |\sum_{l=-\infty}^{k-4}\lambda_lTa_l(x)|>\frac{\lambda}{2}\}|\neq0$, then
$$\lambda\leq 2C_02^{-k(n+\delta+n\delta_2)}\|\chi_k\|_{L^{q'(\cdot)}(\mathbb{R}^{n})}\|f\|_{H\dot{K}^{\alpha,p}_{q(\cdot)}(\mathbb{R}^{n})}.$$
For any given $\lambda>0$, let $k_\lambda$ be the maximal positive integer such that
$$2^{k_\lambda(n+\delta+n\delta_2)}\|\chi_{k_\lambda}\|^{-1}_{L^{q'(\cdot)}(\mathbb{R}^{n})}\leq 2C_0\lambda^{-1}\|f\|_{H\dot{K}^{\alpha,p}_{q(\cdot)}(\mathbb{R}^{n})}.$$
Thus by $\alpha=n\delta_2+\delta$ and Lemma 2.4 we obtain
$$\begin{array}{rl}
\displaystyle V_1&\displaystyle
\leq C\lambda\left\{\sum_{k=-\infty}^{k_\lambda}2^{k\alpha p}\|\chi_{\{x\in A_k: |\sum_{l=-\infty}^{k-4}\lambda_lTa_l(x)|>\frac{\lambda}{2}\}}\|^p_{L^{q(\cdot)}(\mathbb{R}^{n})}\right\}^{1/p}\\
&\displaystyle\leq C\lambda\left\{\sum_{k=-\infty}^{k_\lambda}2^{k\alpha p}\|\chi_k\|_{L^{q(\cdot)}(\mathbb{R}^{n})}^p\right\}^{1/p}\\
&\displaystyle\leq C\lambda\left\{\sum_{k=-\infty}^{k_\lambda}2^{k\alpha p}2^{knp}\|\chi_k\|_{L^{q'(\cdot)}(\mathbb{R}^{n})}^{-p}\right\}^{1/p}\\
&\displaystyle=C\lambda\left\{\sum_{k=-\infty}^{k_\lambda}2^{k(n+\delta+n\delta_2)p}\|\chi_k\|_{L^{q'(\cdot)}(\mathbb{R}^{n})}^{-p}\right\}^{1/p}\\
&\displaystyle\leq C\|f\|_{H\dot{K}^{\alpha,p}_{q(\cdot)}(\mathbb{R}^{n})}.

\end{array}$$

This completes the proof of Theorem 3.3.

\noindent {\bf Remark 3.5}\quad Similar to the method in the proof of Theorem 3.3, we can extend the result of Theorem 3.1 to the case $\alpha=n\delta_2$. Here we omit the details.

If the operator $T$ in Theorem 3.2 is of convolution type, then we can obtain the following stronger conclusion.

\noindent{\bf{Theorem 3.4}}\quad \textit{Let $Tf(x)=\mathrm{p.v.}(k\ast f)(x)$ be a bounded operator on $L^{q(\cdot)}(\mathbb{R}^{n})$ for some $q(\cdot)\in
\mathcal{B}(\mathbb{R}^{n})$, and for some $\delta\in(0,1]$, the kernel $k$ satisfy
$$|k(x-y)-k(x)|\leq C\frac{|y|^\delta}{|x|^{n+\delta}},\quad \mathrm{if}\,\,\,|x|>2|y|.\eqno(3.5)$$
If $n\delta_2\leq\alpha<n\delta_2+\delta$ and $0<p< \infty$, then $T$ can be extended to a bounded operator on $H\dot{K}^{\alpha,p}_{q(\cdot)}(\mathbb{R}^{n})$(or
$HK^{\alpha,p}_{q(\cdot)}(\mathbb{R}^{n})$).}

Using the atomic and molecular theory for the spaces $H\dot{K}^{\alpha,p}_{q(\cdot)}(\mathbb{R}^{n})$ and
$HK^{\alpha,p}_{q(\cdot)}(\mathbb{R}^{n})$, we can prove Theorem 3.4 by a standard procedure, see [15]. Here we omit the details.

For the operator $T$ in Theorem 3.4, we have the following result which is stronger than Theorem 3.3 with the end case $\alpha=n\delta_2+\delta$.

\noindent{\bf{Theorem 3.5}}\quad \textit{Let $T$ be the same as in Theorem 3.4 with $\delta\in(0,1)$. If $\alpha=n\delta_2+\delta$, $q(\cdot)\in
\mathcal{B}(\mathbb{R}^{n})$ and $0<p< \infty$, then $T$ maps continuously $H\dot{K}^{\alpha,p}_{q(\cdot)}(\mathbb{R}^{n})$(or
$HK^{\alpha,p}_{q(\cdot)}(\mathbb{R}^{n})$) into $WH\dot{K}^{\alpha,p}_{q(\cdot)}(\mathbb{R}^{n})$(or
$WHK^{\alpha,p}_{q(\cdot)}(\mathbb{R}^{n})$).}

\noindent{\bf Proof}\quad We only prove the homogeneous case. Let $\alpha=n\delta_2+\delta$ and $f\in
H\dot{K}^{\alpha,p}_{q(\cdot)}(\mathbb{R}^{n})$. By Lemma 2.6 we get
$\displaystyle f=\sum_{k=-\infty}^\infty\lambda_ka_k$, where
$\displaystyle\|f\|_{H\dot{K}^{\alpha,p}_{q(\cdot)}(\mathbb{R}^{n})}\approx
\inf(\sum_{k=-\infty}^\infty|\lambda_k|^{p})^{1/p}$ (the infimum
is taken over above decompositions of $f$), and $a_k$ is a dyadic
central $(\alpha, q(\cdot))$-atom with the support $B_k$. For a fixed $\lambda>0$, we can
write
$$\begin{array}{rl}
&\displaystyle\lambda\left\{\sum_{j=-\infty}^{\infty}2^{j\alpha p}\|\chi_{\{x\in A_j: G_N(Tf)(x)>\lambda\}}\|^p_{L^{q(\cdot)}(\mathbb{R}^{n})}\right\}^{1/p}\\

\end{array}$$
$$\begin{array}{rl}
&\displaystyle\leq C\lambda\left\{\sum_{j=-\infty}^{\infty}2^{j\alpha p}\|\chi_{\{x\in A_j: \sum_{k=-\infty}^{j-4}|\lambda_k|G_N(Ta_k)(x)>\lambda\}}\|^p_{L^{q(\cdot)}(\mathbb{R}^{n})}\right\}^{1/p}\\
&\displaystyle\hspace{3mm}+C\lambda\left\{\sum_{j=-\infty}^{\infty}2^{j\alpha p}\|\chi_{\{x\in A_j: \sum_{k=j-3}^{\infty}|\lambda_k|G_N(Ta_k)(x)>\lambda\}}\|^p_{L^{q(\cdot)}(\mathbb{R}^{n})}\right\}^{1/p}\\
&\displaystyle=E_1+E_2.

\end{array}$$
For $E_2$, by the $L^{q(\cdot)}(\mathbb{R}^{n})$-boundedness of $G_N$ and $T$, we have
$$\begin{array}{rl}
\displaystyle E_2&\displaystyle
\leq C\left\{\sum_{j=-\infty}^{\infty}2^{j\alpha p}\left\|\sum_{k=j-3}^{\infty}|\lambda_k|G_N(Ta_k)\right\|^p_{L^{q(\cdot)}(\mathbb{R}^{n})}\right\}^{1/p}\\
&\displaystyle\leq C\left\{\sum_{j=-\infty}^{\infty}2^{j\alpha p}\left(\sum_{k=j-3}^{\infty}|\lambda_k|^p\|G_N(Ta_k)\|^p_{L^{q(\cdot)}(\mathbb{R}^{n})}\right)\right\}^{1/p}\\
&\displaystyle\leq C\left\{\sum_{j=-\infty}^{\infty}2^{j\alpha p}\left(\sum_{k=j-3}^{\infty}|\lambda_k|^p\|Ta_k\|^p_{L^{q(\cdot)}(\mathbb{R}^{n})}\right)\right\}^{1/p}\\
&\displaystyle\leq C\left\{\sum_{j=-\infty}^{\infty}2^{j\alpha p}\left(\sum_{k=j-3}^{\infty}2^{-k\alpha p}|\lambda_k|^p\right)\right\}^{1/p}\\
&\displaystyle\leq C\left\{\sum_{k=-\infty}^{\infty}|\lambda_k|^p\right\}^{1/p}\leq C\|f\|_{H\dot{K}^{\alpha,p}_{q(\cdot)}(\mathbb{R}^{n})}.

\end{array}$$

Now we consider the term $E_1$. We want to obtain the pointwise estimate for $G_N(Ta_k)(x)$ with $x\in A_j$ and $k\leq j-4$. Note that
$\int_{\mathbb{R}^{n}}Ta_k(x)dx=0$. Let $|x-y|<t$ and write
$$\begin{array}{rl}
\displaystyle|(Ta_k\ast\phi_t)(y)|&\displaystyle\leq\int_{\mathbb{R}^{n}}\left|\int_{B_k}a_k(v)k(w-v)dv\right|t^{-n}\left|\phi\left(\frac{y-w}{t}\right)
-\phi\left(\frac{y}{t}\right)\right|dw\\
&\displaystyle\leq \int_{|w|<2^{k+1}}\left|\int_{B_k}a_k(v)k(w-v)dv\right|t^{-n}\left|\phi\left(\frac{y-w}{t}\right)
-\phi\left(\frac{y}{t}\right)\right|dw\\
&\displaystyle\hspace{3mm}+\int_{2^{k+1}\leq|w|<|x|/2}\left|\int_{B_k}a_k(v)k(w-v)dv\right|t^{-n}\left|\phi\left(\frac{y-w}{t}\right)
-\phi\left(\frac{y}{t}\right)\right|dw\\
&\displaystyle\hspace{3mm}+\int_{|w|\geq|x|/2}\left|\int_{B_k}a_k(v)k(w-v)dv\right|t^{-n}\left|\phi\left(\frac{y-w}{t}\right)
-\phi\left(\frac{y}{t}\right)\right|dw\\
&\displaystyle=:F_{1}+F_{2}+F_{3}.

\end{array}$$
For $F_1$, by Lemma 2.3, the generalized H\"{o}lder inequality and the mean value theorem, we obtain
$$\begin{array}{rl}
\displaystyle F_{1}&\displaystyle\leq C\|Ta_k\|_{L^{q(\cdot)}(\mathbb{R}^{n})}t^{-n}\left\|\left(\phi\left(\frac{y-\cdot}{t}\right)-\phi\left(\frac{y}{t}\right)\right)\chi_{B_{k+1}}\right\|_{L^{q'(\cdot)}(\mathbb{R}^{n})}\\
&\displaystyle\leq C\|a_k\|_{L^{q(\cdot)}(\mathbb{R}^{n})}t^{-n}\left\|\sup_{|\beta|=1}\left|D^\beta\phi\left(\frac{y-\theta\cdot}{t}\right)\right|\frac{|\cdot|}{t}\chi_{B_k}\right\|_{L^{q'(\cdot)}(\mathbb{R}^{n})}\\
&\displaystyle\leq C\|a_k\|_{L^{q(\cdot)}(\mathbb{R}^{n})}\left\|\frac{|\cdot|\chi_{B_k}}{(|x-y|+|y-\theta\cdot|)^{n+1}}\right\|_{L^{q'(\cdot)}(\mathbb{R}^{n})}\\
&\displaystyle\leq C2^{-k\alpha}\frac{1}{|x|^{n+1}}\left\||\cdot|\chi_{B_k}\right\|_{L^{q'(\cdot)}(\mathbb{R}^{n})}\\

\end{array}$$
$$\begin{array}{rl}
&\displaystyle\leq C2^{-k\alpha+k}\frac{1}{|x|^{n+1}}\|\chi_{B_k}\|_{L^{q'(\cdot)}(\mathbb{R}^{n})}\\
&\displaystyle\leq C2^{k\delta-k\alpha}\frac{1}{|x|^{n+\delta}}\|\chi_{B_k}\|_{L^{q'(\cdot)}(\mathbb{R}^{n})}

\end{array}$$
where $0<\theta<1$.
For $F_2$, by (3.5) and the vanishing condition of $a_k$ we have
$$\begin{array}{rl}
\displaystyle F_{2}&\displaystyle=\int_{2^{k+1}\leq|w|<|x|/2}\left|\int_{B_k}a_k(v)(k(w-v)-k(w))dv\right|t^{-n}\left|\phi\left(\frac{y-w}{t}\right)
-\phi\left(\frac{y}{t}\right)\right|dw\\
&\displaystyle\leq\int_{2^{k+1}\leq|w|<|x|/2}\left(\int_{B_k}|a_k(v)||(k(w-v)-k(w)|dv\right)t^{-n}\sup_{|\beta|=1}\left|D^\beta\phi\left(\frac{y-\theta w}{t}\right)\right|\frac{|w|}{t}dw\\
&\displaystyle\leq C\int_{2^{k+1}\leq|w|<|x|/2}\left(\int_{B_k}|a_k(v)||v|^\delta dv\right)\frac{|w|}{|x|^{n+1}}\frac{1}{|w|^{n+\delta}}dw\\
&\displaystyle\leq C2^{k\delta}\|a_k\|_{L^{q(\cdot)}(\mathbb{R}^{n})}\|\chi_{B_k}\|_{L^{q'(\cdot)}(\mathbb{R}^{n})}\frac{1}{|x|^{n+1}}\int_{2^{k+1}\leq|w|<|x|/2}|w|^{1-n-\delta}dw\\
&\displaystyle\leq C2^{k\delta-k\alpha}\|\chi_{B_k}\|_{L^{q'(\cdot)}(\mathbb{R}^{n})}\frac{1}{|x|^{n+1}}|x|^{1-\delta}\\
&\displaystyle\leq C2^{k\delta-k\alpha}\frac{1}{|x|^{n+\delta}}\|\chi_{B_k}\|_{L^{q'(\cdot)}(\mathbb{R}^{n})},

\end{array}$$
where $0<\theta<1$. By (3.5) and the vanishing condition of $a_k$ we can obtain
$$\begin{array}{rl}
\displaystyle F_{3}&\displaystyle\leq\int_{|w|\geq|x|/2}\left|\int_{B_k}a_k(v)(k(w-v)-k(w))dv\right|t^{-n}\left(\left|\phi\left(\frac{y-w}{t}\right)\right|
+\left|\phi\left(\frac{y}{t}\right)\right|\right)dw\\
&\displaystyle\leq C2^{k\delta-k\alpha}\|\chi_{B_k}\|_{L^{q'(\cdot)}(\mathbb{R}^{n})}\int_{|w|\geq|x|/2}|w|^{-n-\delta}t^{-n}\left(\left|\phi\left(\frac{y-w}{t}\right)\right|
+\left|\phi\left(\frac{y}{t}\right)\right|\right)dw\\
&\displaystyle\leq C2^{k\delta-k\alpha}\|\chi_{B_k}\|_{L^{q'(\cdot)}(\mathbb{R}^{n})}\left(|x|^{-n-\delta}+|x|^{-n}\int_{|w|\geq|x|/2}|w|^{-n-\delta}dw\right)\\
&\displaystyle\leq C2^{k\delta-k\alpha}\frac{1}{|x|^{n+\delta}}\|\chi_{B_k}\|_{L^{q'(\cdot)}(\mathbb{R}^{n})},

\end{array}$$
where we have invoked the fact that $t+|y|>|x-y|+|y|>|x|$. Thus, if $x\in A_j$ and $k\leq j-4$, then we have
$$G_N(Ta_k)(x)\leq C2^{k\delta-k\alpha}\frac{1}{|x|^{n+\delta}}\|\chi_{B_k}\|_{L^{q'(\cdot)}(\mathbb{R}^{n})}.$$
Therefore, we obtain
$$\sum_{k=-\infty}^{j-4}|\lambda_k|G_N(Ta_k)(x)\leq C\sum_{k=-\infty}^{j-4}|\lambda_k|2^{k\delta-k\alpha}\frac{1}{|x|^{n+\delta}}\|\chi_{B_k}\|_{L^{q'(\cdot)}(\mathbb{R}^{n})}.$$
If $|\{x\in A_j: \sum_{k=-\infty}^{j-4}|\lambda_k|G_N(Ta_k)(x)>\lambda\}|\neq0$, then by Lemma 2.3 and $\alpha=n\delta_2+\delta$ we have
$$\begin{array}{rl}
\displaystyle \lambda&\displaystyle<C\sum_{k=-\infty}^{j-4}|\lambda_k|2^{k\delta-k\alpha}\frac{1}{|x|^{n+\delta}}\|\chi_{B_k}\|_{L^{q'(\cdot)}(\mathbb{R}^{n})}\\
&\displaystyle\leq C\sum_{k=-\infty}^{j-4}|\lambda_k|2^{k\delta-k\alpha}2^{-j(n+\delta)}2^{(k-j)n\delta_2}\|\chi_{B_j}\|_{L^{q'(\cdot)}(\mathbb{R}^{n})}\\
&\displaystyle\leq C 2^{-j(n+\alpha)}\|\chi_{B_j}\|_{L^{q'(\cdot)}(\mathbb{R}^{n})}\sum_{k=-\infty}^{\infty}|\lambda_k|\\
&\displaystyle\leq C2^{-j(n+\alpha)}\|\chi_{B_j}\|_{L^{q'(\cdot)}(\mathbb{R}^{n})}\|f\|_{H\dot{K}^{\alpha,p}_{q(\cdot)}(\mathbb{R}^{n})}.

\end{array}$$
Let $j_\lambda$ be the maximal integer such that the above inequality holds. Then by Lemma 2.4 we have
$$\begin{array}{rl}
\displaystyle E_1&\displaystyle\leq C\lambda\left(\sum_{j=-\infty}^{j_\lambda}2^{j\alpha p}\|\chi_{j}\|^p_{L^{q(\cdot)}(\mathbb{R}^{n})}\right)^{1/p}\\
&\displaystyle\leq C\lambda\left(\sum_{j=-\infty}^{j_\lambda}2^{j\alpha p}2^{jnp}\|\chi_{B_j}\|^{-p}_{L^{q'(\cdot)}(\mathbb{R}^{n})}\right)^{1/p}\\
&\displaystyle\leq C\|f\|_{H\dot{K}^{\alpha,p}_{q(\cdot)}(\mathbb{R}^{n})}.

\end{array}$$

This finishes the proof of Theorem 3.5.

In general, we have the following theorems. The proofs are similar. Here we omit the details.

\noindent{\bf{Theorem 3.6}}\quad \textit{Let $q(\cdot)\in
\mathcal{B}(\mathbb{R}^{n})$, $\alpha\geq n\delta_2$, $s=[\alpha+n\delta_2]$ and $\varepsilon_0=\alpha+n\delta_2-s$. Suppose $Tf(x)=\mathrm{p.v.}(k\ast f)(x)$ is bounded on $L^{q(\cdot)}(\mathbb{R}^{n})$ and the kernel $k$ satisfies
$$|D^Jk(x-y)-D^Jk(x)|\leq C|y|^\varepsilon|x|^{-n-\varepsilon}$$
for all multi-index $J$ with $|J|=s$, some $\varepsilon>\varepsilon_0$ and $|x|>2|y|$. If $0<p< \infty$,
then $T$ can be extended to a bounded operator on $H\dot{K}^{\alpha,p}_{q(\cdot)}(\mathbb{R}^{n})$(or
$HK^{\alpha,p}_{q(\cdot)}(\mathbb{R}^{n})$).}

\noindent{\bf{Theorem 3.7}}\quad \textit{Let $T$ be the same as in Theorem 3.6 with $\varepsilon=\varepsilon_0$ and $\alpha>n\delta_2$. If $0<p\leq 1$, then $T$ maps continuously $H\dot{K}^{\alpha,p}_{q(\cdot)}(\mathbb{R}^{n})$(or
$HK^{\alpha,p}_{q(\cdot)}(\mathbb{R}^{n})$) into $WH\dot{K}^{\alpha,p}_{q(\cdot)}(\mathbb{R}^{n})$(or
$WHK^{\alpha,p}_{q(\cdot)}(\mathbb{R}^{n})$).}

\medskip

\noindent\bf{Acknowledgements}\quad\rm {The authors are very grateful to the referees
for their valuable comments. This work was supported by National Natural Science Foundation of China (Grant Nos. 11926343, 11926342, 11761026 and 11671397), Shandong Provincial Natural Science Foundation (Grant No. ZR2017MA041) and Project of Shandong Province Higher Educational Science and Technology Program (Grant No. J18KA225).}

\vskip5mm
\centerline{\Large\bf References} \vspace{0.5cm} \def%
\hang{\hangindent\parindent} \def\textindent#1{\indent\llap{#1\enspace}%
\ignorespaces} \def\re{\par\hang\textindent}

\re{[1]} Chen Y, Levin S and Rao M, Variable exponent, linear growth functionals in image restoration, {\it SIAM J. Appl. Math.,} {\bf 66}(2006)1383-1406.

\re{[2]} Cruz-Uribe D, Fiorenza A and Neugebauer C, The
maximal function on variable $L^p$ spaces, {\it Ann. Acad. Sci.
Fenn. Math.,} {\bf 28}(2003)223-238.

\re{[3]} Cruz-Uribe D and Fiorenza A, Variable Lebesgue Spaces: Foundations and Harmonic Analysis(Applied and Numerical
Harmonic Analysis), Springer, Heidelberg, 2013.

\re{[4]} Diening L, Harjulehto P, H\"{a}st\"{o} P and R\r{u}\v{z}i\v{c}ka M, {\it Lebesgue and Sobolev spaces with variable
exponents,} Lecture Notes in Math., vol. {\bf 2017},
Springer, Heidelberg, 2011.

\re{[5]} Ferreira L and P\'{e}rez-L\'{o}pez J, Besov-weak-Herz spaces and global solutions for Navier-Stokes equations, {\it Pacific J. Math.,} {\bf 296}(2018)57-77.

\re{[6]} Harjulehto P, H\"{a}st\"{o} P, L\^{e} U V and Nuortio M, Overview of differential equations with non-standard growth, {\it Nonlinear Anal.,} {\bf 72}(2010)4551-4574.

\re{[7]} Hu G, Lu S and Yang D, The weak Herz spaces. {\it J. Beijing Normal Univ. (Natur. Sci.),} {\bf 33}(1997)27-34.

\re{[8]} Hu G, Lu S and Yang D, The applications of weak Herz spaces. {\it Adv. Math. (China),} {\bf 26}(1997)417-428.

\re{[9]} Izuki M, Boundedness of sublinear operators on Herz spaces
with variable exponent and application to wavelet characterization,
{\it Anal. Math.,} {\bf 36}(2010)33-50.

\re{[10]} Komori Y, Weak type estimates for Calder¨®n-Zygmund operators on Herz spaces at critical indexes,
{\it Math. Nachr.,} {\bf 259}(2003)42-50.

\re{[11]} Kov\'{a}\v{c}ik O and R\'{a}kosn\'{i}k J, On spaces
$L^{p(x)}$ and $W^{k,p(x)}$, {\it Czechoslovak Math. J.,} {\bf
41}(1991)592-618.

\re{[12]} Liu L, The inequalities of commutators on weak Herz spaces,
{\it J. Korean Math. Soc.,} {\bf 39}(2002)899-912.

\re{[13]} R\r{u}\v{z}i\v{c}ka M, Electrorheological fluids: modeling and mathematical theory, Springer, Berlin, 2000.

\re{[14]} Tsutsui Y, The Navier-Stokes equations and weak Herz spaces,
{\it Adv. Differential Equations,} {\bf 16}(2011)1049-1085.

\re{[15]} Wang H and Liu Z, The Herz-type Hardy spaces with variable
exponent and their applications, {\it Taiwanese J. Math.,} {\bf
16}(2012)1363-1389.

\re{[16]} Wang H, Some estimates of intrinsic square functions on
weighted Herz-type Hardy spaces, {\it J. Inequal. Appl.,} {\bf 2015}(2015)62: 22 pages.

\re{[17]} Yan X, Yang D, Yuan W and Zhuo C, Variable weak Hardy spaces and their
applications, {\it J. Funct. Anal.,} {\bf
271}(2016)2822-2887.

\bigskip

\medskip

\noindent\author{Hongbin \uppercase{Wang}}\\
    {School of Mathematics and Statistics, Shandong University of Technology, Zibo, Shandong, 255049, China\\
    E-mail\,$:$ wanghb@sdut.edu.cn}
\medskip

\noindent\author{Zongguang \uppercase{Liu}}\\
    {Department of Mathematics, China University of Mining and Technology(Beijing), Beijing, 100083, China\\
    E-mail\,$:$ liuzg@cumtb.edu.cn}

\end{document}